\documentclass[10pt,oneside,requno]{amsart}

\usepackage{amstext,amsmath,amsfonts,amscd,amsthm,graphicx}
\input xy \xyoption{all}

\theoremstyle{plain}
\newtheorem{theorem}{Theorem}

\newtheorem{proposition}[theorem]{Proposition}

\theoremstyle{definition}
\newtheorem{definition}{Definition}

\theoremstyle{remark}
\newtheorem{remark}{Remark}

\begin{document}

\title[]{Algebra generators and information theory}

\author{Tom\'{a}\v{s} Kopf and Renata Ot\'{a}halov\'{a}}
\address{Mathematical Institute of the Silesian University at Opava\\
Na Rybn\'{\i}\v{c}ku 1 \\
 746 01 Opava, Czech Republic}

\maketitle

\begin{abstract}
A general simplicity problem in category theory is proposed. A
particular example, the simplest choice of generators of an
algebra is specified and illustrated by an example.
\end{abstract}

The idea of physics and of some of mathematics may be seen as the
expansion and compression of information about the real world or, better,
of our understanding of the real world: The understanding has to be
correct and simple.

But what is the correct judgement on simplicity? A vague maximum
principle, Occam's razor is not always well formalized. As a
concept, a measure of simplicity is at hand: The simpler the
thing, the less information is needed to describe it and there is
an accepted measure of information, entropy. This runs however
into at least two problems:

\begin{enumerate}
\item To understand the information, a context is needed. Is the
background information contained in the context also to be
measured? Here we solve this by clearly specifying what structure
is considered background information and what additional
information is measured.

\item A simple-minded evaluation of entropy leads in cases
involving infinite sets or continuous parameters to infinities. In
order to avoid a technically demanding analysis, we count here the
number of continuous parameters to be specified, discounting data
given by finite natural numbers. At the same time we force the
number of parameters to be finite.
\end{enumerate}

The outline of the paper is the following: Section \ref{cat}
measures the information of an object in an abelian category
against a background knowledge category. Section \ref{algebras}
specifies to the case of algebras and sets. An explicit result of
\cite{Otahalova}, Proposition \ref{c} provides an instance of
solving the general problem. A summary is given in the Conclusion.

\section{The simplicity of free resolutions}\label{cat}

The proposal of this section is to measure the information content
of a suitable mathematical structure in an abelian category
according to the following principles:
\begin{itemize}
\item There is a category $A$ and an abelian category $B$ linked
by a forgetful functor $U:B\rightarrow A$ and by a free functor
$F:A\rightarrow B$, the adjoint of $U$. The objects of interest
reside in category $B$ while category $A$ specifies the background
information. \item The functors $U$, $F$ are part of the
background information. Their application is not considered as
connected to the provision of additional information.

\item The information content of an object $b$ of the category $B$
to be measured is the information that needs to be added to $Ub$
in order to understand the structure of b. This is provided in the
form of a free resolution of b. A free resolution \cite{Weibel} is
then given by the following diagram:
\begin{align}
\xymatrix{
 \ar[r]^{{f}_{4}}& F{a}_{3}\ar[r]^{{f}_{3}}& F{a}_{2}\ar[r]^{{f}_{2}}
 & F{a}_{1}\ar[r]^{{f}_{1}}& b\\
 & {a}_{3}\ar[u]^{{i}_{3}}& {a}_{2}\ar[u]^{{i}_{2}}
 & {a}_{1}\ar[u]^{{i}_{1}}&
}\label{freeresolution}
\end{align}
where the upper row is exact. The vertical arrows indicate the functor $F$,
the mappings ${i}_{\bullet}$ should be correctly interpreted as the
inclusions of
generators
${i}_{\bullet}:{a}_{\bullet}\rightarrow UF{a}_{\bullet}$.
\end{itemize}

According to these rules, the information contained in the object $b$
 of the category $B$ with respect to the background knowledge category $A$
is contained in the maps
\begin{align}\label{infomaps}
Uf_k\circ i_k &:a_k\rightarrow UFa_{k-1}&\text{ for $k\in\{ 2, 3, 4,
...\}$}
\end{align}
However, these maps do not need to be specified entirely: There
may be automorphisms of the lower row of the resolution
(\ref{freeresolution}) and respecting the grading of the
resolution that lift to automorphisms of the upper row. Such
automorphisms will be called Bogoliubov automorphisms and
determine the freedom in the maps (\ref{infomaps}) that does not
affect the free resolution (\ref{freeresolution}) and does
therefore not need to be specified. Now, different free
resolutions of $b$ will have different information contents and it
is natural to assign to $b$ the lowest possible information
content attainable by a free resolution.

We may rephrase this by saying that the information content of $b$ is
the information contained in its simplest object $a_1$ of generators
(included in $b$ via the map $Uf_1\circ i_1 : a_1\rightarrow
Ub$) together with its relations as given by the free resolution
(\ref{freeresolution})

Thus, the information to be measured is the information contained in the
maps (\ref{infomaps}) up to Bogoliubov transformations.

That this can be indeed done is shown in an example of the following
section.

\section{Finite free resolutions of algebras and their information
content}\label{algebras}

\subsection{Finite free resolutions}

Let $A$ be the category of vector spaces and $B$ the category of
${C}^{\ast}$-algebras. Let $b$ be a finite dimensional
${C}^{\ast}$-algebra. Then a linear (grading preserving)
transformation transformation in the lower row of
(\ref{freeresolution}) can only be lifted to an algebra
automorphism in the upper row if it maps generators into
generators with the same spectrum. This is obviously not always
the case, for any linear transformation of the lower row of
generators of (\ref{freeresolution}).

It is now possible to enquire into how much information is contained in a
given set of generators or, in other words, in a given free resolution of
the object $b$, an algebra. Unfortunately, the simple idea of just
counting, whatever information
there is and arriving at a quantitative result is hindered by the appearance of
infinities which have to be dealt with properly. To avoid such
technically demanding problems,
we will restrict ourselves to finite dimensional algebras and to finite resolutions
only.

\begin{definition}
A free resolution (\ref{freeresolution}) of an algebra is called here a {\em finite free resolution}
if
\begin{enumerate}
\item there is only a finite number $N$ of  nonzero terms $F{a}_{1}, F{a}_{2}, \ldots , F{a}_{N}$ in the
resolution,
\item the dimensions ${d}_{j}$ of the vector spaces ${a}_{j}$ are all
finite,
\item the maps ${f}_{j}\circ {i}_{j}: {a}_{j}\rightarrow F{a}_{j-1}$, $j=2, 3, 4,\ldots $ send ${a}_{j}$ into a subspace of $F{a}_{j-1}$  of finite degree ${\partial }_{j-1}$ of the filtration of $F{a}_{j-1}$ by word
length.
\end{enumerate}
\end{definition}

\begin{remark}
While the first two conditions of the definition of a finite free resolution are entirely straightforward, the last, third condititon deserves some further justification. The maps ${f}_{i}$ which are an essential part of the free resolution are uniquely determined by the maps ${f}_{j}\circ {i}_{j}: {a}_{j}\rightarrow F{a}_{j-1}$, $j=2, 3, 4,\ldots $ and the information contained in them is to be included into the information contained in the free resolution. In fact, this is  the dominant information content of the resolution. Unless ${a}_{j-1}$ is of dimension ${d}_{j-1}=0$, the free algebra $F{a}_{j-1}$ is necessarily infinite dimensional, spanned as a linear space by (finite) words with letters being elements in ${a}_{j}$ and is canonically a filtered algebra. The linear subspaces ${\mathcal{F}}_{n}F{a}_{j-1}$ of this filtration of  the free algebra $F{a}_{j-1}$ are spanned by words containing at most $n$
letters.

That the range of  ${f}_{j}\circ {i}_{j}: {a}_{j}\rightarrow F{a}_{j-1}$ is inside  is a finite subspace  ${\mathcal{F}}_{{\partial}_{j-1}}F{a}_{j-1}$ is a simple consequence of the finite dimension
${d}_{j}$ of ${a}_{j}$ and of the fact that only finite sums of words
are contained in $F{a}_{j-1}$. Thus the third
requirement is not restrictive but rather serves to define the numbers ${\partial}_{j}$ and to emphasise that infinite series are not included in the
resolution.

A restatement of this may be that all relations (relations of relations, ...) of the resolution ought to be generated from finite polynomial relations (relations of relations,
...).
\end{remark}

\subsection{The information content of a finite free resolution}

To summarize, a finite free resolution is given by the following data:
\begin{itemize}
\item a number $N$ giving the length of the resolution,
\item the numbers ${d}_{1}, \ldots , {d}_{N}$ giving the dimensions of the vector
spaces
${a}_{1}, \ldots , {a}_{N}$,
\item the grading numbers ${\partial}_{1}, \ldots , {\partial}_{N-1}$,
\item the algebra homomorphisms ${f}_{1}, \ldots , {f}_{N}$ determined uniquely by the linear
maps
\begin{align}
{f}_{2}\circ {i}_{2}&: {a}_{2}\rightarrow
{\mathcal{F}}_{{\partial}_{1}}F{a}_{1}\label{firstmap}\\
{f}_{3}\circ {i}_{3}&:
{a}_{3}\rightarrow {\mathcal{F}}_{{\partial}_{2}}F{a}_{2}\\
{f}_{4}\circ {i}_{4}&: {a}_{4}\rightarrow
{\mathcal{F}}_{{\partial}_{3}}F{a}_{3}\\
\intertext{\hskip 5 cm $\vdots$} {f}_{N}\circ {i}_{N}&:
{a}_{N}\rightarrow
{\mathcal{F}}_{{\partial}_{N-1}}F{a}_{N-1}\label{lastmap}
 \end{align}
\end{itemize}

In agreement with \cite{Shannon}, the entropy of a natural number
$N$ will be taken to be
\begin{align}
{S}_{N}&= \ln{N},
\end{align}
and the entropy of all numbers in the data of a finite free resolution is
just the sum of the entropies of  all numbers:
\begin{align}
{S}_{numbers}&= \ln{N} + \sum_{j=1}^{N}{\ln{{d}_{j}}}
+ \sum_{j=1}^{N-1}{\ln{{\partial}_{j}}}.
\end{align}

It remains now to find the entropy of the maps ${f}_{N}\circ
{i}_{N}: {a}_{j}\rightarrow
{\mathcal{F}}_{{\partial}_{N-1}}F{a}_{N-1}$ modulo Bogoljubov
automorphisms. That is the essential part of the entropy
calculation, since specifying these maps requires to determine real
parameters. These contain an infinite number of digits and carry
therefore an infinite entropy against which one typically can neglect
${S}_{numbers}$. Thus the minimalization of required information boils
down to finding the set of maps (\ref{firstmap})-(\ref{lastmap})
containing up to  Bogoliubov automorphisms the least number of continuous
parameter specifications.

This problem has been solved in \cite{Otahalova} for the case of $b$
being the algebra $M_2(\mathbb{C})$ of $2\times 2$-matrices with complex
numbers in its entries. The solution is found to be the Clifford
resolution of $M_2(\mathbb{C})$:

\begin{definition}
Let $M_{2^m}(\mathbb{C})$ be the algebra of ${2^m}\times
{2^m}$-matrices. Its {\bf Clifford resolution}  is the following finite
free resolution:
\begin{align}
\xymatrix{
 0 \ar[r]^{}& F{a}_{2}\ar[r]^{f_2}
 & F{a}_{1}\ar[r]^{\gamma}& M_{2^m}(\mathbb{C})\\
 & V_{m(2m+1)}\ar[u]^{{i}_{2}}
 & V_{2m}\ar[u]^{{i}_{1}}&
}\label{freeresolution}
\end{align}
Where $V_{2m}$ is the $2m$-dimensional complex vector space of
generators with a basis $e_1$, $e_2$, ..., $e_{2m}$ and where
$V_{m(2m+1)}$ is the $m(2m+1)$-dimensional complex vector space of
relations with a basis $r_{1\leq 1}$, $r_{1\leq 2}$, ...,
$r_{2m\leq 2m}$ mapped by $f_2\circ i_2$ onto the relations
\begin{align}
f_2\circ i_2 r_{k\leq l} &= e_k e_l + e_l e_k - {\delta}_{kl}
&\text{ for $k,l \in \{ 1,2,...,2m\}$ and $k\leq l$}
\end{align}
\end{definition}

\begin{proposition}\label{c}
The finite free resolution of the matrix algebra $M_2(\mathbb{C})$
with the least number of parameters to be specified is given by
the Clifford resolution.
\end{proposition}
For the convenience of the reader we provide the following
\begin{proof}
Since the algebra $M_2(\mathbb{C})$ is noncommutative, its
subspace of generators has to be at least 2-dimensional.

The word length of relations has to be at least 2, again to avoid
commutativity.

The algebra $M_2(\mathbb{C})$ is 4-dimensional, with the
2-dimensional set of generators and the unit being linearly
independent, leaving thus one dimension for one linearly
independent linear combination of four words of length 2. Thus at
least 3 relations are required.

The Clifford resolution allows for a 1-dimensional Bogoliubov
automorphism. This can be checked to be the maximal dimension for
any possible Bogoliubov automorphism.
\end{proof}

\section{Conclusion}
The general proposal of Section \ref{cat} to measure information
content of the object of a category through its simplest free
resolution against a category of background knowledge was
specified in Section \ref{algebras} for the case of
$C^{\ast}$-algebras and their linear subspaces of generators. In
the simplest case of $2\times 2$-matrices, the simplest finite
free resolution is known \cite{Otahalova}. That result adds as an
example on substance to the general framework explained here and
poses at the same time the challenge of more such examples.

Note that the simplest finite free resolution is of interest not
just as a tool to measure information but also for providing a
suitable set of generators well adapted to the algebra in
question.

The case of $n\times n$-matrices will be discussed in further work
of one of the authors (R. Ot\'{a}halov\'{a}).

\section{Acknowledgements}
We thank Michal Marvan for comments on an earlier draft. This work
was supported by grant MSM:J10/98:192400002.


\end{document}